# RELATIONS BETWEEN THE STRONG GLOBAL DIMENSION, COMPLEXES OF FIZED SIZE AND DERIVED CATEGORY

YOHNY CALDERÓN-HENAO, FELIPE GALLEGO-OLAYA,
AND HERNÁN GIRALDO

ABSTRACT. Let $\mathbb{Z}$ be the integer numbers, $\mathbb{K}$ an algebraically closed field, $\Lambda$ a finite dimensional $\mathbb{K}$-algebra, mod$\Lambda$ the category of finitely generated right modules, proj$\Lambda$ the full subcategory of mod$\Lambda$ consisting of all projective $\Lambda$-modules, and $C_n(\text{proj}\Lambda)$ the bounded complexes of projective $\Lambda$-modules of fixed size for any integer $n \geq 2$. We find an algorithm to calculate the strong global dimension of $\Lambda$, when $\Lambda$ is a finite strong global dimension and derived discrete, using the Auslander-Reiten quivers of the categories $C_n(\text{proj}\Lambda)$. Also, we show the relation between the Auslander-Reiten quiver of the bounded derived category $D^b(\Lambda)$ and the Auslander-Reiten quiver of $C_{\eta+1}(\text{proj}\Lambda)$, where $\eta = $ s.gl.dim$(\Lambda)$ (strong global dimension of $\Lambda$).

## INTRODUCTION

Let $\mathbb{Z}$ be the integer numbers, $\mathbb{K}$ an algebraically closed field, $\Lambda$ a finite dimensional $\mathbb{K}$-algebra, and mod$\Lambda$ the category of finitely generated right $\Lambda$-modules. The theory of almost split sequences (Auslander-Reiten sequence) was extended with success to the study of the bounded derived category of mod$\Lambda$, denoted by $D^b(\Lambda)$. The notion of almost split sequences gave origin to the notion of almost split triangles (Auslander-Reiten triangles). In his fundamental book [13], D. Happel proved that the derived category of a finite global dimension algebra has almost split triangles, he was very successful in describing the Auslander-Reiten quiver in the case of a hereditary algebra.

Let proj$\Lambda$ be the full subcategory of mod$\Lambda$ consisting of all projective $\Lambda$-modules. The category of all complexes of finitely generated projective right $\Lambda$-modules of fixed size, denoted by $C_n(\text{proj}\Lambda)$ for any integer $n \geq 2$, was introduced and studied by R. Bautista, M.J. Souto Salorio, and R. Zuazua in [4]. In [10], C. Chaio, I. Pratti, and M.J. Souto Salorio continued with the study of these categories, and they showed that the Auslander-Reiten quiver of $C_n(\text{proj}\Lambda)$ can be constructed by the well-known knitting algorithm used to build the Auslander-Reiten quiver mod$\Lambda$. Furthermore, in [3], R. Bautista studied the representation type of $C_n(\text{proj}\Lambda)$ for all positive integer $n \geq 2$. Recently in [9], C. Chaio, A.G. Chaio, and I. Pratti, proved that if the strong global dimension of $\Lambda$ is finite ($\eta = $ s.gl.dim$(\Lambda)$ with $\eta \in \mathbb{N}$), then, $C_{\eta+1}(\text{proj}\Lambda)$ is of finite representation type if and only if $C_n(\text{proj}\Lambda)$ is of finite representation type for each $n \geq 2$.

At the moment, there is not a general method to calculate the strong global dimension of $\Lambda$, but some works have been done in this direction, see [2] and [9]. In the first one, E.R. Alvares, P. Le Meur, and E.N. Marcos used the structure of the derived category to characterize the strong global dimension of $\Lambda$, when $\Lambda$ is piecewise hereditary algebra. In the second one, C. Chaio, A.G. Chaio, and I. Pratti determined the strong global dimension of some picewise hereditary finite dimensional algebras taking into account their ordinary quivers with relations. In this work, we find an algorithm to calculate the strong global dimension of $\Lambda$, when $\Lambda$ is finite strong global dimension and derived discrete. For this, we use the techniques given in [10] for constructing the Auslander-Reiten quiver of the category $C_n(\text{proj}\Lambda)$ for $n \geq 2$, the Theorem 2.18, and the Algorithm 2.19. Also, we proved that if the strong global dimension of $\Lambda$ is finite then $C_{\eta+1}(\text{proj}\Lambda)$ is of tame representation type if and only if $C_n(\text{proj}\Lambda)$ is of tame representation type for all $n \geq 1$ (see Theorem 3.5).

Furthermore, we show that for a finite strong global dimensional algebra $\Lambda$, if we have irreducible morphisms and $\mathcal{E}_{\eta+1}$-almost split sequences in $C_{\eta+1}(\text{proj}\Lambda)$ if and only if we have "almost the same" irreducible morphisms and $\mathcal{E}_{\eta+i}$-almost split sequences in $C_{\eta+i}(\text{proj}\Lambda)$ for $i \geq 2$ (see Theorem 2.16) and also irreducible morphisms and Auslander-Reiten triangle in $D^b(\Lambda)$ (see Theorem 3.1). Finally, we get that the Auslander-Reiten quiver of $D^b(\Lambda)$ is $\mathbb{Z}$ copies of a special subquiver of the Auslander-Reiten of $C_{\eta+1}(\text{proj}\Lambda)$ (see Definition 3.3 and Theorem 3.4). Additionally we note that in [8], C. Chaio, A.G. Chaio, and I. Pratti, using the knitting technique to build the Auslander-Reiten quiver of the category $C_n(\text{proj}\Lambda)$, described the complexes that belong to the mouth of non-homogeneous tubes in the Auslander-Reiten quiver of their bounded derived category, whenever this algebras are either derived equivalent to hereditary algebras of type $\widetilde{A}_n$ or $\widetilde{D}_n$.

This paper is organized as follows: In Section 1 we recall some notations and preliminary results. In Section 2 (the main ideas and more details of this section can be found in [11]), we show how the indecomposable complexes, the irreducible morphisms, and the almost split sequences in $C_{\eta+i}(\text{proj}\Lambda)$ for $i \geq 1$ are when $\Lambda$ is a finite strong global dimension. As a result of our research, we obtain an algorithm to calculate the strong global dimension of a finite dimension $\mathbb{K}$-algebra. In Section 3, we show the relationship between $C_{\eta+1}(\text{proj}\Lambda)$ and $D^b(\Lambda)$, we also provide two examples for appling the algorithm and calculting the strong global dimensional. Furthermore, we get the Auslander-Reiten quiver of the bounded derived category $D^b(\Lambda)$ from the Auslander-Reiten quiver of $C_{\eta+1}(\text{proj}\Lambda)$.

1. Preliminary results

First we will set some notations and definitions. Let $\mathbb{Z}$ be the integer numbers, $\mathbb{K}$ an algebraically closed field, and $\Lambda$ a finite dimensional $\mathbb{K}$-algebra. By $\text{mod}\Lambda$ and $\text{proj}\Lambda$, we denote the category of finitely generated right $\Lambda$-modules and the full subcategory of finitely generated projective right $\Lambda$-modules, respectively. Given a quiver $\Gamma = (\Gamma_0, \Gamma_1)$, where $\Gamma_0$ is the set of vertices and $\Gamma_1$ is the set of arrows, $\alpha\colon i \to j$ will denote an arrow from the vertex $i$ to the vertex $j$. The notation $P_i$, $I_i$, $S_i$ is used for the indecomposable projective (resp. injective, simple) associated to the vertex $i \in \Gamma_0$.

$\mathcal{A}$ will always denote an additive category which is Krull-Schmidt and $\text{Hom}_\mathcal{A}(X, Y)$ the set of the morphism from $X$ to $Y$, where $X$ and $Y$ are objects in $\mathcal{A}$, when no confusion arises we will use the notation $\text{Hom}(X, Y)$. A morphism $f\colon X \to Y$ in $\text{Hom}_\mathcal{A}(X, Y)$ is said to be a section (resp. retraction) if and only if there is a morphism $h\colon Y \to X$ such that $hf = 1_X$ (resp. $fh = 1_Y$). When one of these conditions holds, $f$ is said to be a split morphism. A morphism $f$ of $\mathcal{A}$ is called **irreducible** when it is not split but, for any factorization $f = hg$, $h$ is a retraction or $g$ is a section. A morphism $f\colon E \longrightarrow M$ in $\mathcal{A}$ is called **right almost split** if $f$ is not a retraction and if $g\colon X \longrightarrow M$ is not a retraction, there is a $s\colon X \longrightarrow E$ with $fs = g$. A morphism $f\colon E \longrightarrow M$ in $\mathcal{A}$ is called **minimal right** if $fu = f$ with $u \in \text{End}_\mathcal{A}(E)$ implies $u \in \text{Aut}_\mathcal{A}(E)$. A morphism $f\colon E \longrightarrow M$ in $\mathcal{A}$ is called **minimal right almost split** if $f$ is right almost split and $f$ is minimal right. Dually, a **minimal left almost split morphism** is defined. A pair $(i, d)$ of composable morphism $X \xrightarrow{i} Y \xrightarrow{d} Z$ in $\mathcal{A}$ is said to be **exact** if $i$ is a kernel of $d$ and $d$ is a cokernel of $i$. Let $\mathcal{E}$ be class of exact pairs closed under isomorphisms. The pairs in $\mathcal{E}$ are called **conflactions**. A conflaction in $\mathcal{A}$, $M \xrightarrow{f} E \xrightarrow{g} N$ is said to be **almost split** if $f$ is a minimal left almost split morphism and $g$ is a minimal right almost split morphism.

We recall that a complex $X = (X^i, d_X^i)_{i \in \mathbb{Z}}$ is a family of objects $X^i$ in $\mathcal{A}$ and morphisms $d_X^i \colon X^i \to X^{i+1}$ in $\mathcal{A}$, $i \in \mathbb{Z}$, such that $d_X^{i+1} d_X^i = 0$ for all $i \in \mathbb{Z}$. If $X$ and $Y$ are complexes, a morphism of complexes $f \colon X \to Y$ is given by a family of morphisms $f^i \colon X^i \longrightarrow Y^i$, $i \in \mathbb{Z}$, such that $d_Y^i f^i = f^{i+1} d_X^i$ for all $i \in \mathbb{Z}$. We denote by $\mathrm{C}(\mathcal{A})$ the category whose objects are the complexes over $\mathcal{A}$, the morphisms are the morphism between two complexes, each $X^i$ is called a **cell** and each $f^i$ is called **component** of $f$ for each $i \in \mathbb{Z}$ and $\mathrm{C}^{-,b}(\mathcal{A})$ the full subcategory of $\mathrm{C}(\mathcal{A})$ where the complexes are bounded above and bounded cohomology. We denote by $\mathrm{C}^b(\mathcal{A})$ the full subcategory of $\mathrm{C}(\mathcal{A})$ given by the bounded complexes, by $\mathrm{K}(\mathcal{A})$, $\mathrm{K}^b(\mathcal{A})$ and $\mathrm{K}^{-,b}(\mathcal{A})$ the homotopy categories of $\mathrm{C}(\mathcal{A})$, $\mathrm{C}^b(\mathcal{A})$ and $\mathrm{C}^{-,b}(\mathcal{A})$ respectively, and by $\mathrm{D}^b(\mathcal{A})$ the bounded derived category of $\mathcal{A}$. The bounded derived category of mod-$\Lambda$ is denoted by $\mathrm{D}^b(\Lambda)$. If $X$ is a complex and $p \in \mathbb{Z}$, we have a new complex $X[p]$ given by the shifting or translation functor.

From [14], we recall the following. If $X \in \mathrm{K}^b(\mathcal{A})$ is a complex we may consider a preimage $\bar{X}$ of $X$ in $\mathrm{C}^b(\mathcal{A})$ without indecomposable projective direct summands. $\bar{X}$ is uniquely determined by $X$ up to isomorphism of bounded complexes in $\mathrm{C}^b(\mathcal{A})$. Thus the following is well-defined. If $0 \neq X \in \mathrm{K}^b(\mathcal{A})$ there exists $r \leq s$ such that $\bar{X}^r \neq 0 \neq \bar{X}^s$ and $\bar{X}^i = 0$ for $i < r$ and $i > s$. We will always identify $\bar{X}$ with $X$, we can also identify the complex $X$ in $\mathrm{C}_{(s-r+1)}(\mathrm{proj}\Lambda)$ in the following way: $X^1 = X^r$, $X^i = X^{r+(i-1)}$ for $i \in \{2, \ldots, s-r\}$ and $X^{s-r+1} = X^s$. The length of $X$ is defined as $\ell(X) = s - r$ and the strong global dimensional of $\Lambda$ is defined by

$$\mathrm{s.gl.dim}(\Lambda) \;:=\; \sup\{\ell(X) | X \in \mathrm{K}^b(\mathrm{proj}\Lambda) \text{ indecomposable}\}.$$

From [13], we say that a finite $\mathbb{K}$-algebra $\Lambda$ is to be piecewise hereditary, if there exists an hereditary abelian category $\mathcal{H}$ such that the bounded derived categories $\mathrm{D}^b(\Lambda)$ and $\mathrm{D}^b(\mathcal{H})$ are equivalent as triangulated categories.

The following theorem is proved in [[14], Theorem 3.2].

**Theorem 1.1.** *Let $\Lambda$ be a finite dimensional $\mathbb{K}$-algebra. Then $\Lambda$ is a piecewise hereditary if and only if $\mathrm{s.gl.dim}(\Lambda) < \infty$.*

For a positive integer $n \geq 2$, the category $\mathrm{C}_n(\mathrm{proj}\Lambda)$ consists of the complex $X$, such that $X^j = 0$ if $j \notin \{1, \ldots, n\}$. A complex $X$ in $\mathrm{C}_n(\mathrm{proj}\Lambda)$ with $n \geq 2$ has the following shape

$$X \;:\; \cdots \;\to\; 0 \;\to\; X^1 \;\to\; X^2 \;\to\; \cdots \;\to\; X^{n-1} \;\to\; X^n \;\to\; 0 \;\to\; \cdots$$

when no confusion arises this complex has a fixed sixe, we denoted it by

$$X \;:\; X^1 \;\to\; X^2 \;\to\; \cdots \;\to\; X^{n-1} \;\to\; X^n.$$

And with this notation, we call $X^1$ the first cell, $X^2$ the second cell, ..., $X^n$ the last cell.

In [4], the authors proved that in $\mathrm{C}_n(\mathrm{proj}\Lambda)$ the class of composable morphisms $X \xrightarrow{f} Y \xrightarrow{g} Z$ such that for $j \in \{1, \ldots, n\}$ the sequence $0 \longrightarrow X^j \xrightarrow{f^j} Y^j \xrightarrow{g^j} Z^j \longrightarrow 0$ is split is a conflation in $\mathrm{C}_n(\mathrm{proj}\Lambda)$, and they showed the existence of the $\mathcal{E}_n$-almost split sequences or almost split sequences in $\mathrm{C}_n(\mathrm{proj}\Lambda)$. Also they defined the following complexes in $\mathrm{C}_n(\mathrm{proj}\Lambda)$. Given an indecomposable $P \in \mathrm{proj}\Lambda$

$$\mathrm{S}(P) \colon P \longrightarrow \cdots \longrightarrow 0 \longrightarrow 0, \quad \mathrm{T}(P) \colon 0 \longrightarrow \cdots \longrightarrow 0 \longrightarrow P, \text{ and}$$

$$\mathrm{J}_k(P) \colon 0 \longrightarrow \cdots \longrightarrow 0 \longrightarrow P \xrightarrow{1} P \longrightarrow 0 \longrightarrow \cdots \longrightarrow 0,$$

where $k \in \{1,\ldots, n-1\}$. Let $Y, Z$ be complexes in $C_n(\text{proj}\Lambda)$. A morphism $v\colon Y \to Z$ in $C_n(\text{proj}\Lambda)$ is called **proper epimorphism** if there is an exact sequence $0 \to X \to Y \to Z \to 0$ in $C_n(\text{mod}\Lambda)$. Dually, a **proper monomorphism** is defined. A complex $P$ in $C_n(\text{proj}\Lambda)$ is called $\mathcal{E}_n$**-projective** if for all proper epimorphisms $v\colon Y \to Z$ and any morphism $f\colon P \to Z$ in $C_n(\text{proj}\Lambda)$ there is $g\colon P \to Y$ such that $f = gv$. Dually, $\mathcal{E}_n$**-injective** is defined. They proved the following proposition [[4], Corollary 3.8 and Corollary 3.9]

**Proposition 1.2.** *The indecomposable complexes $\mathcal{E}_n$-**projective** (resp $\mathcal{E}_n$-**injective**) in $C_n(\text{proj}\Lambda)$ are the complexex $J_k(P)$ and $T(P)$ (resp $J_k(P)$ and $S(P)$) for $k \in \{1,\ldots, n-1\}$ with $P$ an indecomposable projective $\Lambda$-module.*

Let $X$ be a complex in $C_{n+1}(\text{proj}\Lambda)$. We denote by $_*[X]$ the complex in $C_n(\text{proj}\Lambda)$ given for $_*[X] : X^2 \to \cdots \to X^n \to X^{n+1}$, similarly, we denote by $[X]_*$ the complex in $C_n(\text{proj}\Lambda)$ given for $[X]_* : X^1 \to X^2 \to \cdots \to X^n$. For a morphism $f\colon X \to Y$ in $C_{n+1}(\text{proj}\Lambda)$ we have

1.
$$\begin{array}{ccccccc}
_*[X] : & X^2 & \to & \cdots & \to & X^n & \to & X^{n+1} \\
& \downarrow _*[f] & \downarrow f^2 & & & \downarrow f^n & & \downarrow f^{n+1} \\
_*[Y] : & Y^2 & \to & \cdots & \to & Y^n & \to & Y^{n+1}.
\end{array}$$

2.
$$\begin{array}{ccccccc}
[X]_* : & X^1 & \to & X^2 & \to & \cdots & \to & X^n \\
& \downarrow [f]_* & \downarrow f^1 & & \downarrow f^2 & & & \downarrow f^n \\
[Y]_* : & Y^1 & \to & Y^2 & \to & \cdots & \to & Y^n.
\end{array}$$

The functor $_*[-]$ is called **left abrupt truncation** and the functor $[-]_*$ is called **right abrupt truncation**.

**Lemma 1.3.** *Let $\Lambda$ be a finite dimensional $\mathbb{K}$-algebra, $X, Y$ be complexes in $C_n(\text{proj}\Lambda)$, and $f : X \longrightarrow Y$ is an irreducible morphism in $C_n(\text{proj}\Lambda)$. The following conditions hold.*

i. *If $X^1 = 0 = Y^1$, then $_*[f] : {_*[X]} \longrightarrow {_*[Y]}$ is an irreducible morphism in $C_{n-1}(\text{proj}\Lambda)$.*

ii. *If $X^n = 0 = Y^n$, then $[f]_* : [X]_* \longrightarrow [Y]_*$ is an irreducible morphism in $C_{n-1}(\text{proj}\Lambda)$.*

**Proof.** This fact follows from the definition of functor. □

Let $X$ be a complex in $C_n(\text{proj}\Lambda)$. We denoted by $I_R(X)$ the complex in $C_{n+1}(\text{proj}\Lambda)$ given for $X^1 \to \cdots \to X^n \to 0$; similarly, we denote by $_L I(X)$ the complex in $C_{n+1}(\text{proj}\Lambda)$ given for $0 \to X^1 \to \cdots \to X^n$. For a morphism $f\colon X \to Y$ in $C_n(\text{proj}\Lambda)$ we have

1.
$$\begin{array}{ccccccc}
I_R(X) : & X^1 & \to & X^2 & \to & \cdots & \to & X^n & \to & 0 \\
& \downarrow I_R(f) & \downarrow f^1 & & \downarrow f^2 & & & \downarrow f^n & & \downarrow \\
I_R(Y) : & Y^1 & \to & Y^2 & \to & \cdots & \to & Y^n & \to & 0.
\end{array}$$

2.
$$\begin{array}{ccccccc}
_L I(X) : & 0 & \to & X^1 & \to & \cdots & \to & X^{n-1} & \to & X^n \\
& \downarrow _L I(X) & \downarrow & & \downarrow f^1 & & & \downarrow f^{n-1} & & \downarrow f^n \\
_L I(X) : & 0 & \to & Y^1 & \to & \cdots & \to & Y^{n-1} & \to & Y^n.
\end{array}$$

The functor $_L\mathrm{I}(-)$ is called **left inclusion** and the functor $\mathrm{I}_R(-)$ is called **right inclusion**. We note that $_L\mathrm{I}(-)$ is the shifting or traslation functor [1] and $\mathrm{I}_R(-)$ is the natural inclusion from $\mathrm{C}_n(\mathrm{proj}\Lambda)$ to $\mathrm{C}_{n+1}(\mathrm{proj}\Lambda)$.

From [8], we get the following definition. Let $X$ be an indecomposable complex in $\mathrm{C}_n(\mathrm{proj}\Lambda)$, $X: X^1 \xrightarrow{d_X^1} X^2 \to \cdots X^{n-1} \xrightarrow{d_X^{n-1}} X^n$. We say that $X$ can be **extended to the left** if there is a projective $\Lambda$-module $X^0$ and a non-zero morphism $\delta : X^0 \longrightarrow X^1$ such that $d_X^1 \delta = 0$. Dually, we say that $X$ can be **extended to the right**. They showed the following result, which has been modified to above notation.

**Proposition 1.4.** *[[8], Lemma 4.3] Let $X$ and $Y$ be complexes in $\mathrm{C}_n(\mathrm{proj}\Lambda)$. Let $f : X \longrightarrow Y$ be an irreducible morphism in $\mathrm{C}_n(\mathrm{proj}\Lambda)$ with $n \geq 2$. Then, the following conditions hold.*

  i. *If $X$ is an indecomposable and it can not be extended to the left then $_L\mathrm{I}(f) : {_L\mathrm{I}(X)} \longrightarrow {_L\mathrm{I}(Y)}$ is an irreducible morphism in $\mathrm{C}_{n+1}(\mathrm{proj}\Lambda)$.*

  ii. *If $Y$ is an indecomposable and it can not be extended to the right then $\mathrm{I}_R(f) : \mathrm{I}_R(X) \longrightarrow \mathrm{I}_R(Y)$ is an irreducible morphism in $\mathrm{C}_{n+1}(\mathrm{proj}\Lambda)$.*

The following proposition wich is the result of joining Corollary 2 and the Proposition 3 proved by H. Giraldo and H. Merklen in [12].

**Proposition 1.5.** *Let $f: X \to Y$ be an irreducible morphism in $\mathrm{C}_n(\mathrm{proj}\Lambda)$. One of the next conditios hold:*

  i. *$f^n$ is a section for all $n \in \mathbb{Z}$.*

  ii. *$f^n$ is a retraction for all $n \in \mathbb{Z}$.*

  iii. *there is an $i_0 \in \mathbb{Z}$ such that $f^{i_0}$ is irreducible, the morphisms $f^j$ are sections for all $j > i_0$ and the morphism $f^j$ are retractions for all $j < i_0$.*

## 2. Main results

For the remainder of this section, $\Lambda$ denoted a piecewise hereditary algebra, that is the same by Theorem 1.1, $\Lambda$ is finite strong global dimension algebra, which is denoted by $\eta = \mathrm{s.gl.dim}(\Lambda)$, with $\eta \geq 1$.

**Remark 2.1.**

  i. If $X \in \mathrm{K}^b(\mathrm{proj}\Lambda)$ is an indecomposable complex with $\ell(X) = \eta$, then X is (up to traslation) neither $\mathcal{E}_{(\eta+1)}$-indecomposable projective nor $\mathcal{E}_{(\eta+1)}$-indecomposable injective.

  In fact, suppose $X$ is either $\mathcal{E}_{(\eta+1)}$-indecomposable projective or $\mathcal{E}_{(\eta+1)}$-indecomposable injective, then by Proposition 1.2 $X$ is some of the following complexes in $\mathrm{C}_{\eta+1}(\mathrm{proj}\Lambda)$

  $$(1) \quad \mathrm{T}(P), \quad (2) \quad \mathrm{S}(P), \quad \text{or} \quad (3) \quad \mathrm{J}_k(P) \text{ con } k \in \{1, \ldots, \eta\},$$

  for $P$ an indecomposable projective $\Lambda$-module. To the complex $\mathrm{T}(P)$ and $\mathrm{S}(P)$

  $$\ell(\mathrm{T}(P)) = 0 = \ell(\mathrm{S}(P))$$

  and for the complex $\mathrm{J}_k(P)$ with $k \in \{1, \ldots, \eta\}$ is zero in $\mathrm{K}^b(\mathrm{proj}\Lambda)$. Therefore in any case $\ell(X) = 0$.

  ii. Let $X$ be indecomposable complex in $\mathrm{K}^b(\mathrm{proj}\Lambda)$ non-zero with $\ell(X) = \eta$ therefore there are $r, s \in \mathbb{Z}$ with $r < s$ such that $X^r \neq 0, X^s \neq 0$ y $X^i = 0$ to $i > s$ and $i < r$. Then $X^j \neq 0$ for $j \in \{r, r+1, \ldots, s\}$.

   iii. Let $X$ be an indecomposable complex in $K^b(\text{proj}\Lambda)$ with $\ell(X) = \eta$, then up to traslation $X \in C_{\eta+i}(\text{proj}\Lambda)$ for $i \geq 1$ and $X \notin C_\eta(\text{proj}\Lambda)$.

**Lemma 2.2.** *Every indecomposable complex in $C_{\eta+2}(\text{proj}\Lambda)$, is a complex such that either the cell one or cell $\eta+2$ are zero.*

**Proof.** We suppose $X$ is an indecomposable complex in $C_{\eta+2}(\text{proj}\Lambda)$, then either $X$ is zero in $K^b(\text{proj}\Lambda)$ or $X$ is non-zero in $K^b(\text{proj}\Lambda)$.

If $X$ is zero in $K^b(\text{proj}\Lambda)$, then $X = J_k(P)$ with $P$ being an indecomposable projecive $\Lambda$-module for $k \in \{1, \ldots, \eta+1\}$. We noted that $X$ satisfies the condition.

On the other hand, if $X$ is non-zero in $K^b(\text{proj}\Lambda)$, we suppose $X^1 \neq 0$ and $X^{\eta+2} \neq 0$, then $\ell(X) = \eta + 1$ and in this way we can say s.gl.dim$(\Lambda) \geq \eta + 1$, contradicting the fact that $\eta = $ s.gl.dim$(\Lambda)$. $\square$

**Corollary 2.3.** *Every indecomposable complex in $C_{\eta+i}(\text{proj}\Lambda)$ is a complex such that either the cell one or cell $\eta+i$ are zero for $i \geq 2$.*

**Proposition 2.4.** *In $C_{\eta+2}(\text{proj}\Lambda)$. If $f : X \longrightarrow Y$ is an irreducible morphism and $X, Y$ are indecomposable complexes, then $X^1 = 0 = Y^1$ or $X^{\eta+2} = 0 = Y^{\eta+2}$, that is, $f$ has the following shapes:*

1.
$$\begin{array}{ccccccc} X & : & X^1 & \longrightarrow & \cdots & \longrightarrow & X^{\eta+1} & \longrightarrow & 0 \\ \downarrow f & & \downarrow & & & & \downarrow & & \downarrow \\ Y & & Y^1 & \longrightarrow & \cdots & \longrightarrow & Y^{\eta+1} & \longrightarrow & 0 \end{array}, \text{ or}$$

2.
$$\begin{array}{ccccccc} X & : & 0 & \longrightarrow & X^2 & \longrightarrow & \cdots & \longrightarrow & X^{\eta+2} \\ \downarrow f & & \downarrow & & \downarrow & & & & \downarrow \\ Y & : & 0 & \longrightarrow & Y^2 & \longrightarrow & \cdots & \longrightarrow & Y^{\eta+2} \end{array}.$$

**Proof.** Since $\eta = $ s.gl.dim$(\Lambda)$ and $X, Y$ are indecomposable complexes in $C_{\eta+2}(\text{proj}\Lambda)$, by Lemma 2.2 we can consider all possible forms for $X$ and $Y$ to $f : X \longrightarrow Y$.

We show that $f : X \longrightarrow Y$ can only be of the shapes above and $f$ does not have the next two shapes

$$\begin{array}{ccccccc} X & : & 0 & \longrightarrow & \tilde{X}^1 & \longrightarrow & \cdots & \longrightarrow & \tilde{X}^\eta & \longrightarrow & \tilde{X}^{\eta+1} \\ \downarrow f & & \downarrow & & \downarrow & & & & \downarrow & & \downarrow \\ Y & : & Y^1 & \longrightarrow & Y^2 & \longrightarrow & \cdots & \longrightarrow & Y^{\eta+1} & \longrightarrow & 0 \end{array} \text{ or}$$

$$\begin{array}{ccccccc} X & : & X^1 & \longrightarrow & X^2 & \longrightarrow & \cdots & \longrightarrow & X^{\eta+1} & \longrightarrow & 0 \\ \downarrow f & & \downarrow & & \downarrow & & & & \downarrow & & \downarrow \\ Y & : & 0 & \longrightarrow & \tilde{Y}^1 & \longrightarrow & \cdots & \longrightarrow & \tilde{Y}^\eta & \longrightarrow & \tilde{Y}^{\eta+1} \end{array}.$$

If we have the form:

$$\begin{array}{ccccccc} X & : & 0 & \longrightarrow & \tilde{X}^1 & \longrightarrow & \cdots & \longrightarrow & \tilde{X}^\eta & \longrightarrow & \tilde{X}^{\eta+1} \\ \downarrow f & & \downarrow & & \downarrow & & & & \downarrow & & \downarrow \\ Y & : & Y^1 & \longrightarrow & Y^2 & \longrightarrow & \cdots & \longrightarrow & Y^{\eta+1} & \longrightarrow & 0. \end{array}$$

  a) If $\tilde{X}^{\eta+1} = 0$, then this case is the case 1.

b) If $Y^1 = 0$, then this case is the case 2.

Now if $\tilde{X}^{\eta+1} \neq 0$ and $Y^1 \neq 0$, by Proposition 1.5, how as $f$ is an irreducible morphism in $C_{\eta+2}(\text{proj}\Lambda)$ then one of these conditions holds:

i. For each $i \in \{1, \ldots, \eta+2\}$, the morphisms $f^i$ are sections in $\text{proj}\Lambda$.

ii. For each $i \in \{1, \ldots, \eta+2\}$, the morphisms $f^i$ are retractions in $\text{proj}\Lambda$.

iii. There is an $i \in \{1, \ldots, \eta+2\}$ such that $f^i$ is irreducible in $\text{proj}\Lambda$, the morphisms $f^j$ are sections for all $j > i$ and the morphisms $f^j$ are retractions for all $j < i$.

But $f^{\eta+2} : \tilde{X}^{\eta+1} \longrightarrow 0$ is not a section in $\text{proj}\Lambda$ and the morphism $f^1 : 0 \longrightarrow Y^1$ is not retraction in $\text{proj}\Lambda$. Therefore i, ii y iii are not satisfied, which is a contradiction.

If we have the form:

$$\begin{array}{ccccccccc}
X & : & X^1 & \longrightarrow & X^2 & \longrightarrow & \cdots & \longrightarrow & X^{\eta+1} & \longrightarrow & 0 \\
& & \downarrow f & & \downarrow & & & & \downarrow & & \downarrow \\
Y & : & 0 & \longrightarrow & \tilde{Y}^1 & \longrightarrow & \cdots & \longrightarrow & \tilde{Y}^\eta & \longrightarrow & \tilde{Y}^{\eta+1}
\end{array}$$

a) If $X^1 = 0$, this case is the case 2.

b) If $\tilde{Y}^{\eta+1} = 0$, this case is the case 1.

Now if $X^1 \neq 0$ and $\tilde{Y}^{\eta+1} \neq 0$.

By hypothesis $f$ is an irreducible morphism in $C_{\eta+2}(\text{proj}\Lambda)$ and $X$ and $Y$ indecomposable complexes and we note that $X$ and $Y$ are not the form $J_k(P)$ with $P$ being an indecomposable projective $\Lambda$-module. Using Theorem 5 stated in [12] we have that $X$ and $Y$ are minimal projective complexes and by Theorem 6 stated in [12] $f$ is an irreducible morphism in $K^b(\text{proj}\Lambda)$, therefore by the first proposition stated in [1], section 5.1, the cone of an irreducible morphism in $K^b(\text{proj}\Lambda)$ is indecomposable, that is, $C_f \colon X^1 \longrightarrow X^2 \longrightarrow \cdots \longrightarrow \tilde{Y}^\eta \longrightarrow \tilde{Y}^{\eta+1}$ is indecomposable in $K^b(\text{proj}\Lambda)$. Since $X^1 \neq 0$ and $\tilde{Y}^{\eta+1} \neq 0$, then $\ell(C_f) = \eta + 1$, contradicting the fact that $\eta = \text{s.gl.dim}(\Lambda)$. □

**Corollary 2.5.** *In $C_{\eta+i}(\text{proj}\Lambda)$ with $i \geq 2$. Let $X, Y$ be indecomposable complexes and $f \colon X \longrightarrow Y$ be an irreducible morphism.*

a) *If $\ell(X) = l - j = \eta$, then $Y^{j-r} = Y^{l+r} = 0$ for all $r \geq 1$.*

b) *If $\ell(Y) = l - j = \eta$, then $X^{j-r} = X^{l+r} = 0$ for all $r \geq 1$.*

That is, with $\eta = l - j$, so if $\ell(X) = \eta$ or $\ell(Y) = \eta$ then $f$ has the following form.

$$\begin{array}{ccccccccccc}
X \colon & 0 & \longrightarrow & \cdots & \longrightarrow & 0 & \longrightarrow & X^j & \longrightarrow & \cdots & \longrightarrow & X^l & \longrightarrow & 0 & \longrightarrow & \cdots & \longrightarrow & 0 \\
& \downarrow f & & & & & & \downarrow f^j & & & & \downarrow f^l & & & & & \\
Y \colon & 0 & \longrightarrow & \cdots & \longrightarrow & 0 & \longrightarrow & Y^j & \longrightarrow & \cdots & \longrightarrow & X^l & \longrightarrow & 0 & \longrightarrow & \cdots & \longrightarrow & 0.
\end{array}$$

**Proof.** We get a) and b) in a way similar to what was done in the proof of Proposition 2.4. □

**Corollary 2.6.** *Let $X, Y$ be indecomposable complexes in $C_{\eta+2}(\mathrm{proj}\Lambda)$. If $f : X \longrightarrow Y$ is an irreducible morphism in $C_{\eta+2}(\mathrm{proj}\Lambda)$, then $_*[f]$ or $[f]_*$ is an irreducible morphism in $C_{\eta+1}(\mathrm{proj}\Lambda)$.*

**Proof.** By Proposition 2.4 $X^1 = 0 = Y^1$ or $X^{\eta+2} = 0 = Y^{\eta+2}$.

If $X^1 = 0 = Y^1$, then by Lemma 1.3 $_*[f]$ is an irreducible morphism in $C_{\eta+1}(\mathrm{proj}\Lambda)$.

If $X^{\eta+2} = 0 = Y^{\eta+2}$, then by Lemma 1.3 $[f]_*$ is an irreducible morphism in $C_{\eta+1}(\mathrm{proj}\Lambda)$.

□

**Remark 2.7.** From Corollary 2.6 and Corollary 2.3, we note that if $f: X \to Y$ is an irreducible morphism in $C_{\eta+i}(\mathrm{proj}\Lambda)$ for $i \geq 2$, with $X$ and $Y$ indecomposable complexes, we have "almost the same $f$", an irreducible morphism in $C_{\eta+1}(\mathrm{proj}\Lambda)$ up to applying the left or right abrupt truncations $(i-1)$-times. For this reason if $f$ is an irreducible morphism in $C_{\eta+i}(\mathrm{proj}\Lambda)$ with $i \geq 2$, we say that $f$ is an irreducible morphism in $C_{\eta+1}(\mathrm{proj}\Lambda)$.

**Lemma 2.8.** *Let $X, Z$ be indecomposables complexes in $C_{\eta+2}(\mathrm{proj}\Lambda)$ and $\delta : X \longrightarrow Y \longrightarrow Z$ an almost split sequence in $C_{\eta+2}(\mathrm{proj}\Lambda)$, then*

  i. *If $X^1 = 0 = Y^1$ and $Z^1 = 0$, then $_*[\delta] : {}_*[X] \longrightarrow {}_*[Y] \longrightarrow {}_*[Z]$ is an almost split sequence in $C_{\eta+1}(\mathrm{proj}\Lambda)$.*

  ii. *If $X^{\eta+2} = 0 = Y^{\eta+2}$ and $Z^{\eta+2} = 0$, then $[\delta]_* : [X]_* \longrightarrow [Y]_* \longrightarrow [Z]_*$ is an almost split sequence in $C_{\eta+1}(\mathrm{proj}\Lambda)$.*

**Proof.** Let $\delta : X \xrightarrow{f} Y \xrightarrow{g} Z$ be an almost split sequence in $C_{\eta+2}(\mathrm{proj}\Lambda)$.

  i. We assume $X^1 = 0 = Y^1$ and $Z^1 = 0$, we know that $f$ is a minimal left almost split morphism in $C_{\eta+2}(\mathrm{proj}\Lambda)$, in particular $f$ is an irreducible morphism in $C_{\eta+2}(\mathrm{proj}\Lambda)$. By the Lemma 1.3 $_*[f] : {}_*[X] \longrightarrow {}_*[Y]$ is an irreducible morphism in $C_{\eta+1}(\mathrm{proj}\Lambda)$.

  We will prove that $_*[f] : {}_*[X] \longrightarrow {}_*[Y]$ and $_*[g] : {}_*[Y] \longrightarrow {}_*[Z]$ are a minimal left and right almost split morphisms in $C_{\eta+1}(\mathrm{proj}\Lambda)$ respectively.

  1. Since $_*[f] : {}_*[X] \longrightarrow {}_*[Y]$ is an irreducible morphism in $C_{\eta+1}(\mathrm{proj}\Lambda)$, then $_*[f]$ is not a section in $C_{\eta+1}(\mathrm{proj}\Lambda)$.

  2. Let $u : {}_*[X] \longrightarrow U$ be a morphism in $C_{\eta+1}(\mathrm{proj}\Lambda)$ such that $u$ is not a section in $C_{\eta+1}(\mathrm{proj}\Lambda)$, we will prove that there is $v : {}_*[Y] \longrightarrow U$ a morphism in $C_{\eta+1}(\mathrm{proj}\Lambda)$ such that $v \circ {}_*[f] = u$. We considered the complex $_LI(U)$ and the morphism $_LI(u)$ in $C_{\eta+2}(\mathrm{proj}\Lambda)$. We note that $_LI(u)$ is not a section in $C_{\eta+2}(\mathrm{proj}\Lambda)$, since $u$ is not a section in $C_{\eta+1}(\mathrm{proj}\Lambda)$. As $f$ is a minimal left almost split morphism in $C_{\eta+2}(\mathrm{proj}\Lambda)$, there is $h : Y \longrightarrow {}_LI(U)$ a morphism in $C_{\eta+2}(\mathrm{proj}\Lambda)$, such that $h \circ f = {}_LI(u)$. If we considered the following morphism $v : {}_*[Y] \longrightarrow U$ in $C_{\eta+1}(\mathrm{proj}\Lambda)$ which is defined by $v := (h^2, \ldots, h^{\eta+2})$, then $v \circ {}_*[f] = u$.

  3. Let $r \in \mathrm{End}({}_*[Y])$ be such that $r \circ {}_*[f] = {}_*[f]$, so we have $_LI(r) \in \mathrm{End}(Y)$ and $_LI(r) f = f$, and since $f$ is minimal in $C_{\eta+2}(\mathrm{proj}\Lambda)$, then $_LI(r) \in \mathrm{Aut}(Y)$ and therefore $r \in \mathrm{Aut}({}_*[Y])$.

  From 1, 2 and 3, $_*[f]$ is a minimal left almost split morphism in $C_{\eta+1}(\mathrm{proj}\Lambda)$.

In a similar way we can prove that $_*[g] : {_*[Y]} \longrightarrow {_*[Z]}$ is minimal right almost split morphism in $\mathrm{C}_{\eta+1}(\mathrm{proj}\Lambda)$.

This proves that $_*[\delta]\colon {_*[X]} \longrightarrow {_*[Y]} \longrightarrow {_*[Z]}$ is an almost split sequence in $\mathrm{C}_{\eta+1}(\mathrm{proj}\Lambda)$.

ii. The proof is similar to i.

□

**Proposition 2.9.** *Let $X, Z$ be indecomposables complexes in $\mathrm{C}_{\eta+2}(\mathrm{proj}\Lambda)$ and $\delta\colon X \xrightarrow{f} Y \xrightarrow{g} Z$ an almost split sequence in $\mathrm{C}_{\eta+2}(\mathrm{proj}\Lambda)$, then $_*[\delta]$ or $[\delta]_*$ is an almost split secuence in $\mathrm{C}_{\eta+1}(\mathrm{proj}\Lambda)$.*

**Proof.** Since $X \xrightarrow{f} Y$ is an irreducible morphism in $\mathrm{C}_{\eta+2}(\mathrm{proj}\Lambda)$ by Proposition 2.4 we have that $X^1 = 0 = Y^1$ or $X^{\eta+2} = 0 = Y^{\eta+2}$. If $X^1 = 0 = Y^1$, then $Z^1 = 0$ since $0 \longrightarrow X^1 \longrightarrow Y^1 \longrightarrow Z^1 \longrightarrow 0$ is split. By Lemma 2.8 $_*[\delta]$ is an almost split secuence in $\mathrm{C}_{\eta+1}(\mathrm{proj}\Lambda)$. If $X^{\eta+2} = 0 = Y^{\eta+2}$, similarly $Z^{\eta+2} = 0$ and by Lemma 2.8 $[\delta]_*$ is an almost split secuence in $\mathrm{C}_{\eta+1}(\mathrm{proj}\Lambda)$.

□

**Remark 2.10.** From Proposition 2.9 and Remark 2.7, we note that if $\delta$ is an almost split sequence in $\mathrm{C}_{\eta+i}(\mathrm{proj}\Lambda)$ with $i \geq 2$, we have "almost the same $\delta$" an almost split sequence in $\mathrm{C}_{\eta+1}(\mathrm{proj}\Lambda)$ up to applying the left or right abrupt truncations $(i-1)$-times. For this reason if $\delta$ is an almost split sequence in $\mathrm{C}_{\eta+i}(\mathrm{proj}\Lambda)$ we say that $\delta$ is an almost split sequence in $\mathrm{C}_{\eta+1}(\mathrm{proj}\Lambda)$.

**Lemma 2.11.** *Let $X \in \mathrm{C}_{\eta+1}(\mathrm{proj}\Lambda)$ be an indecomposable complex.*

*i. If $X$ can be extended to the left, then $X^{\eta+1} = 0$.*

*ii. If $X$ can be extended to the right, then $X^1 = 0$.*

**Proof.** Let $X \in \mathrm{C}_{\eta+1}(\mathrm{proj}\Lambda)$ be an indecomposable complex. For both cases i and ii, we note that $X \neq \mathrm{J}_k(P)$ with $k \in \{1, \ldots, \eta\}$ and $P$ being an indecomposable projective $\Lambda$-module.

i. We assume that $X^{\eta+1} \neq 0$ and $X$ can be extended to the left, so there is $X^0 \in \mathrm{proj}\Lambda$ and non-zero morphism $d_X^0\colon X^0 \to X^1$ such that $d_X^1 d_X^0 = 0$. We consider the morphism $f$ in $\mathrm{K}^b(\mathrm{proj}\Lambda)$

$$\begin{array}{ccccccccccc}
X' : & \cdots & \longrightarrow & X^0 & \longrightarrow & 0 & \longrightarrow & \cdots & \longrightarrow & 0 & \longrightarrow & 0 & \longrightarrow & \cdots \\
& & & \downarrow f & & \downarrow d_X^0 & & \downarrow & & & & \downarrow & & \downarrow & & \\
X : & \cdots & \longrightarrow & X^1 & \xrightarrow{d_X^1} & X^2 & \longrightarrow & \cdots & \longrightarrow & X^\eta & \xrightarrow{d_X^\eta} & X^{\eta+1} & \longrightarrow & \cdots
\end{array}$$

We note that $f$ is non-zero and not invertible, since $X$ is indecomposable. How $\mathrm{Hom}_{\mathrm{K}(\mathrm{proj}\Lambda)}(X, X'[1]) = 0$ by [[14], Corollary 1.4], then $\mathrm{C}_f$ is indecomposable, where $\mathrm{C}_f = X^0 \to X^1 \to \cdots \to X^\eta \to X^{\eta+1}$, then $\ell(\mathrm{C}_f) = \eta + 1$. But this is a contradiction, because $\eta = \mathrm{s.gl.dim}(\Lambda)$.

ii. Similar to i.

$\square$

**Proposition 2.12.** *Let $X, Y \in \mathrm{C}_{\eta+1}(\mathrm{proj}\Lambda)$ be indecomposable complexes and $f: X \to Y$ an irreducible morphism in $\mathrm{C}_{\eta+1}(\mathrm{proj}\Lambda)$, then $_L\mathrm{I}(f)$ or $\mathrm{I}_R(f)$ is an irreducible morphism in $\mathrm{C}_{\eta+2}(\mathrm{proj}\Lambda)$.*

**Proof.**

i. If $X$ can not be extended to the left, by Proposition 1.4, $_L\mathrm{I}(f)$ it is an irreducible morphism in $\mathrm{C}_{\eta+2}(\mathrm{proj}\Lambda)$.

ii. If $X$ can be extended to the left, there is a projective $\Lambda$-module $X^0$ and a non-zero morphism $d_X^0 : X^0 \longrightarrow X^1$ such that $d_X^1 d_X^0 = 0$. We note that $d_X^0$ is not epimorphism, if $d_X^0$ were epimorphism then $X$ would be a decomposable complex. We consider the following composable morphisms $X^0 \xrightarrow{d_X^0} X^1 \xrightarrow{f^1} Y^1$.

   A) If $f^1 d_X^0 \neq 0$, how $d_Y^1(f^1 d_X^0) = f^2(d_X^1 d_X^0) = f^2 \circ 0 = 0$, then $Y$ can be extended to the left, and by Lemma 2.11, $Y^{\eta+1} = 0$. Therefore $Y$ can not be extended to the right and by Proposition 1.4, $\mathrm{I}_R(f)$ is an irreducible morphism in $\mathrm{C}_{\eta+2}(\mathrm{proj}\Lambda)$.

   B) If $f^1 d_X^0 = 0$, we have $\mathrm{im} d_X^0 \subseteq \ker f^1$ and since $d_X^0$ is a non-zero morphism, then $Y^1 \neq 0$ and by Lemma 2.11, thus $Y$ can not be extended to the right. By Proposition 1.4, $\mathrm{I}_R(f)$ is an irreducible morphism in $\mathrm{C}_{\eta+2}(\mathrm{proj}\Lambda)$.

$\square$

**Remark 2.13.**

i. From Proposition 2.12 and Lemma 2.11, we note that if $f$ is an irreducible morphism in $\mathrm{C}_{\eta+1}(\mathrm{proj}\Lambda)$, we have "almost the same $f$", an irreducible morphism in $\mathrm{C}_{\eta+i}(\mathrm{proj}\Lambda)$ with $i \geq 2$ up to applying the left or right inclusions $(i-1)$-times, for this reason if $f$ is an irreducible morphism in $\mathrm{C}_{\eta+1}(\mathrm{proj}\Lambda)$ we say that $f$ is an irreducible morphism in $\mathrm{C}_{\eta+i}(\mathrm{proj}\Lambda)$ with $i \geq 2$.

ii. If $X, Y \in \mathrm{C}_{\eta+1}(\mathrm{proj}\Lambda)$ are indecomposable complexes and $f: X \to Y$ is an irreducible morphism in $\mathrm{C}_{\eta+1}(\mathrm{proj}\Lambda)$. From the proof of Proposition 2.12, if $X$ can be extended to the left then $Y$ can not be extended to the right.

**Proposition 2.14.** *Let $X, Z$ be indecomposables complexes in $\mathrm{C}_{\eta+1}(\mathrm{proj}\Lambda)$ and $\delta: X \xrightarrow{f} Y \xrightarrow{g} Z$ is an almost split sequence in $\mathrm{C}_{\eta+1}(\mathrm{proj}\Lambda)$, then $_L\mathrm{I}(\delta)$ or $\mathrm{I}_R(\delta)$ is an almost split sequence in $\mathrm{C}_{\eta+2}(\mathrm{proj}\Lambda)$.*

**Proof.** Let $\delta: X \xrightarrow{f} Y \xrightarrow{g} Z$ be an almost split sequence in $\mathrm{C}_{\eta+1}(\mathrm{proj}\Lambda)$

i. If $X$ can not be extended to the left. By Proposition 1.4 each $_L\mathrm{I}(f):_L\mathrm{I}(X) \to _L\mathrm{I}(Y)$ is an irreducible morphism in $\mathrm{C}_{\eta+2}(\mathrm{proj}\Lambda)$. We will show that $_L\mathrm{I}(f)$ is a minimal left almost split morphims. In fact, let $v:_L\mathrm{I}(X) \to A$ be morphism in $\mathrm{C}_{\eta+2}(\mathrm{proj}\Lambda)$ such that $v$ is not a section, we note that $_*[v]$ is not a section in $\mathrm{C}_{\eta+1}(\mathrm{proj}\Lambda)$, since $X$ can not be extended to the left. How $f$ is a minimal left almost split morphism in $\mathrm{C}_{\eta+1}(\mathrm{proj}\Lambda)$ there is $\tilde{v}: Y \to \tilde{A}$ such that $\tilde{v}f =_* [v]$. If we define $h := (0, \tilde{v}^1, \ldots, \tilde{v}^{\eta+2})$ then $h \circ_L \mathrm{I}(f) = v$. This shows that $_L\mathrm{I}(f)$ is left almost split morphism. Using left abrupt truncation functor shows that $_L\mathrm{I}(f)$ is a minimal morphism, therefore $_L\mathrm{I}(f)$ is a minimal left almost split morphism in $\mathrm{C}_{\eta+2}(\mathrm{proj}\Lambda)$ and so $_L\mathrm{I}(\delta)$ is an almost split sequence in $\mathrm{C}_{\eta+2}(\mathrm{proj}\Lambda)$.

ii. If $X$ can be extended to the left. Let $Y = \oplus_{i=1}^{m} Y_i$, where $Y_i$ is indecomposable for all $i \in \{1, \ldots, m\}$. So $f_i : X \to Y_i$ is an irreducible morphism with $f = \begin{bmatrix} f_1 \\ \vdots \\ f_m \end{bmatrix}$, and by Remark 2.13 ii, $Y_i$ can not be extended to the right for all $i \in \{1, \ldots, m\}$, therefore by Proposition 1.4, $I_R(f_i) : I_R(X) \to I_R(Y_i)$ is an irreducible morphism in $C_{\eta+2}(\text{proj}\Lambda)$ for all $i \in \{1, \ldots, m\}$, so $I_R(f) = \begin{bmatrix} I_R(f_1) \\ \vdots \\ I_R(f_m) \end{bmatrix}$ is an irreducible morphism in $C_{\eta+2}(\text{proj}\Lambda)$. Similarly to was done in item i, it is shown that $I_R(f)$ is a minimal left almost split morphism in $C_{\eta+2}(\text{proj}\Lambda)$, so $I_R(\delta)$ is an almost split secuence in $C_{\eta+2}(\text{proj}\Lambda)$.

□

**Remark 2.15.** From Proposition 2.14 and Remark 2.13 i., we note that if $\delta$ is an almost split sequence in $C_{\eta+1}(\text{proj}\Lambda)$, we have "almost the same $\delta$", an almost split sequence in $C_{\eta+i}(\text{proj}\Lambda)$ with $i \geq 2$ up to applying the left or right inclusions $(i-1)$-times. For this reason if $\delta$ is an almost split sequence in $C_{\eta+1}(\text{proj}\Lambda)$ we say that $\delta$ is an almost split sequence in $C_{\eta+i}(\text{proj}\Lambda)$ with $i \geq 2$.

**Theorem 2.16.** Let $X, Y$ be indecomposables complexes in $C_{\eta+1}(\text{proj}\Lambda)$. Then, the following conditions hold:

i. $f : X \longrightarrow Y$ is an irreducible morphism in $C_{\eta+2}(\text{proj}\Lambda)$ if and only if $f : X \longrightarrow Y$ is an irreducible morphism in $C_{\eta+1}(\text{proj}\Lambda)$.

ii. $\delta : X \longrightarrow B \longrightarrow Y$ is an almost split sequence in $C_{\eta+2}(\text{proj}\Lambda)$ if and only $\delta : X \longrightarrow B \longrightarrow Y$ is an almost split sequence in $C_{\eta+1}(\text{proj}\Lambda)$.

**Proof.**

i. ($\Longrightarrow$) Corollary 2.6.
   ($\Longleftarrow$) Proposition 2.12.

ii. ($\Longrightarrow$) Proposition 2.9.
   ($\Longleftarrow$) Proposition 2.14.

□

As an immediate consequence of the above result and by Remark 2.7, Remark 2.10, Remark 2.13 i. and Remark 2.15, we can state the following corollary of the form.

**Corollary 2.17.** Let $X, Y$ be indecomposable complexes in $C_{\eta+1}(\text{proj}\Lambda)$. $f : X \longrightarrow Y$ is an irreducible morphism in $C_{\eta+1}(\text{proj}\Lambda)$ and $\delta : X \longrightarrow B \longrightarrow Y$ is an almost split sequence in $C_{\eta+1}(\text{proj}\Lambda)$.

i. $f : X \longrightarrow Y$ is an irreducible morphism in $C_{\eta+1}(\text{proj}\Lambda)$ if and only if $f$ is an irreducible morphism in $C_{\eta+i}(\text{proj}\Lambda)$ with $i \geq 2$.

ii. $\delta : X \longrightarrow B \longrightarrow Y$ is an almost split sequence in $C_{\eta+1}(\text{proj}\Lambda)$ if and only if $\delta$ is an almost split sequence in $C_{\eta+i}(\text{proj}\Lambda)$ with $i \geq 2$.

The following result give us the new method to calculate the strong global dimension for some finite dimensional $\mathbb{K}$-algebras.

**Theorem 2.18.** *If $\Lambda$ is a finite dimensional $\mathbb{K}$-algebra, then* $\mathrm{s.gl.dim}(\Lambda) = m_0 - 2$ *where*

$$m_0 := \min\{n \geq 2 \,|\, \forall X \in \mathrm{C}_n(\mathrm{proj}\Lambda) \text{ indecomposable}, X^1 = 0 \text{ or } X^n = 0\}.$$

**Proof.** We note the following

i. If $\mathrm{s.gl.dim}(\Lambda) = \infty$, the result is clear.
ii. If $\mathrm{s.gl.dim}(\Lambda) = \eta$, then there is $X$ indecomposable complex in $\mathrm{C}_{\eta+1}(\mathrm{proj}\Lambda)$ such that $\ell(X) = \eta$. The complex $X$ satisfies that $X^1 \neq 0$ and $X^{\eta+1} \neq 0$, therefore

$$\eta + 1 \;<\; m_0.$$

And finally by Lemma 2.2, we have $m_0 \leq \eta + 2$. Therefore $\eta = m_0 - 2$.

$\square$

Following [16], a derived category $\mathrm{D}^b(\Lambda)$ is said to be derived discrete if for every vector $\boldsymbol{n} = (n_i)_{i \in \mathbb{Z}}$ of natural numbers there are only finitely many isomorphism classes of indecomposable objects in $\mathrm{D}^b(\Lambda)$ of homology dimension vector $\boldsymbol{n}$.

If $\Lambda$ is derived discrete by [[3], Theorem 2.3 (a)], then $\mathrm{C}_n(\mathrm{proj}\Lambda)$ is of finite representation type for all $n \geq 1$. Furthermore, if $\Lambda$ is piecewise hereditary by definition it exists an hereditary abelian category $\mathcal{H}$ such that $\mathrm{D}^b(\Lambda) \cong \mathrm{D}^b(\mathcal{H})$ as triangulated categories. So, from [14] and [16], if $\Lambda$ is piecewise hereditary and derived discrete, then $\mathcal{H}$ is given by a path algebra of Dynkin type, and we can calculate the strong global dimension of $\Lambda$ explicited using the following algorithm.

**Algorithm 2.19.** Let $\Lambda$ be a piecewise hereditary and derived discrete. In [10], the authors show how to build the Auslander-Reiten quiver of the category $\mathrm{C}_n(\mathrm{proj}\Lambda)$ for all $n \geq 1$.

Step 1. Build the Auslander-Reiten quiver of $\mathrm{C}_2(\mathrm{proj}\Lambda)$.

a) If all indecomposable complexes satisfy that the first and the last cells are zero, then by Theorem 2.18 we have $\mathrm{s.gl.dim}(\Lambda) = 2 - 2 = 0$.
b) Otherwise, go to step 2.

Step 2. Build the Auslander-Reiten quiver of $\mathrm{C}_3(\mathrm{proj}\Lambda)$ and do the a) or b) of the Step 1.

Step 3. Continue until you obtain the case (a) of the previous step. The algorithm stops because $\Lambda$ is a finite strong global dimension by Theorem 1.1.

## 3. Applications for the derived category and examples

First we start showing the relation between $\mathrm{C}_{\eta+1}(\mathrm{proj}\Lambda)$ and the bounded derived category $\mathrm{D}^b(\Lambda)$, when $\eta = \mathrm{s.gl.dim}(\Lambda)$. Also, in this section, we will always use the well-know fact $\mathrm{K}^{-,b}(\mathrm{proj}\Lambda) \cong \mathrm{D}^b(\Lambda)$.

**Theorem 3.1.** *Let $\Lambda$ be a piecewise hereditary, with $\eta = \mathrm{s.gl.dim}(\Lambda)$. Let be $X$ and $Z$ indecomposable complexes in $\mathrm{C}_{\eta+1}(\mathrm{proj}\Lambda)$ where $X, Z$ can not be neither a $\mathcal{E}_{\eta+1}$-projective nor a $\mathcal{E}_{\eta+1}$-injective.*

i. *If $X$ can not be extended to the left and $Z$ can not be extended to the right, $f: X \to Z$ is an irreducible morphism in $\mathrm{C}_{\eta+1}(\mathrm{proj}\Lambda)$ if and only if $f: X \to Z$ is an irreducible morphism in $\mathrm{D}^b(\Lambda)$.*

    ii. Let $Y = \oplus_{i=1}^{m} Y_i$ with $Y_i$ an indecomposable complex in $C_{\eta+1}(\text{proj}\Lambda)$ such that can not be neither a $\mathcal{E}_{\eta+1}$-projective nor a $\mathcal{E}_{\eta+1}$-injective for all $i \in \{1, \ldots, m\}$, $X$ can not be extended to the left, $Y_i$ neither can not be extended to the left nor to the right for all $i \in \{1, \ldots, m\}$, and $Z$ can not be extended to the right. Then, $\delta$: $X \xrightarrow{f} Y \oplus P \xrightarrow{g} Z$ is an almost split sequence in $C_{\eta+1}(\text{proj}\Lambda)$ where $P$ is both $\mathcal{E}_{\eta+1}$-projective and $\mathcal{E}_{\eta+1}$-injective indecomposable complex if and only if there is $\delta'$: $X \xrightarrow{f'} Y \xrightarrow{g'} Z$ is Auslander-Reiten triangle in $D^b(\Lambda)$.

**Proof.**

i. See [[8], Theorem 4.4].

ii. Since $\delta \colon X \xrightarrow{f} Y \oplus P \xrightarrow{g} Z$ is an almost split sequence in $C_{\eta+1}(\text{proj}\Lambda)$, then $f$ is a minimal left almost split morphism in $C_{\eta+1}(\text{proj}\Lambda)$, therefore $f_i$ is an irreducible morphism in $C_{\eta+1}(\text{proj}\Lambda)$ for all $i \in \{1, \ldots, m\}$ where $f_m \colon X \to P$. Then by i. $f_i \colon X \to Y_i$ with $i \in \{1, \ldots, m-1\}$ is an irreducible morphism in $D^b(\Lambda)$. We note that $f' \colon X \longrightarrow Y$ with $f' = (f_1, \ldots, f_{m-1})^t$ is a minimal left almost split morphism in a way similar to what was done in the proof of Proposition 2.14 i. so $\delta' \colon X \xrightarrow{f'} Y' \xrightarrow{g'} Z$ is Auslander-Reiten triangle in $D^b(\Lambda)$.

The converse is follows from Theorem 2.7 in [7], which is also true in the $D^b(\Lambda)$ (see [6]).

□

**Remark 3.2.** The result [[8], Theorem 4.4] holds for $\mathbb{K}$-algebras of finite dimension.

We will define an special subquiver of Auslander-Reiten quiver of $C_{\eta+1}(\text{proj}\Lambda)$.

**Definition 3.3.** *Let $\Lambda$ be a piecewise hereditary algebra and $\eta = \text{s.gl.dim}(\Lambda)$. Let $\Gamma_{C_{\eta+1}(\text{proj}\Lambda)}$ be the Auslander-Reiten quiver of $C_{\eta+1}(\text{proj}\Lambda)$, we define $\bar{\Gamma}_{C_{\eta+1}(\text{proj}\Lambda)}$ by the quiver obtained from $\Gamma_{C_{\eta+1}(\text{proj}\Lambda)}$ in the following way, in each connected component of $\Gamma_{C_{\eta+1}(\text{proj}\Lambda)}$ we remove the indecomposable complexes such that are both $\mathcal{E}_{\eta+1}$-projective and $\mathcal{E}_{\eta+1}$-injective and afterwards we take the connected component given by the indecomposables complexes $X$ such that $X$ can not be extended niether to the left nor to the right.*

**Theorem 3.4.** *With the above notation. Let $\Lambda$ be a piecewise hereditary algebra and $\eta = \text{s.gl.dim}(\Lambda)$ the Auslander-Reiten quiver of $D^b(\Lambda)$ is $\mathbb{Z}$ copies of $\bar{\Gamma}_{C_{\eta+1}(\text{proj}\Lambda)}$.*

**Proof.** The result follows from Theorem 3.1, Corollary 2.5, the definition of piecewise hereditary algebra, and the well known fact: the Auslander-Reiten quiver of $D^b(\mathcal{H})$ is $\mathbb{Z}$ copies of Auslander-Reiten quiver of $\text{mod}\mathcal{H}$, where $\mathcal{H}$ is an hereditary algebra. The latter is, $\Gamma_{D^b(\Lambda)} = \mathbb{Z}\Gamma_{\text{mod}(\mathcal{H})}$, where $\Gamma_{D^b(\Lambda)}$ and $\Gamma_{\text{mod}(\mathcal{H})}$ are the Auslander-Reiten quivers of $D^b(\Lambda)$ and $\text{mod}(\mathcal{H})$ respectively [[13], Chapter I.5.5. Corollary].

□

**Theorem 3.5.** *Let $\Lambda$ be a piecewise hereditary with $\eta = \text{s.gl.dim}(\Lambda)$. $C_{\eta+1}(\text{proj}\Lambda)$ is of tame representation type if and only if $C_n(\text{proj}\Lambda)$ is of tame representation type for all $n \geq 1$.*

**Proof.** Suppose that $C_{\eta+1}(\text{proj}\Lambda)$ is of tame representation type.

i. If $n \leq \eta + 1$, by [[3], Theorem 1.1] $C_n(\text{proj}\Lambda)$ is of tame representation type or wild representation type. If $C_n(\text{proj}\Lambda)$ is of wild representation type, how $C_n(\text{proj}\Lambda) \subseteq C_{\eta+1}(\text{proj}\Lambda)$, then $C_{\eta+1}(\text{proj}\Lambda)$ is of wild representation type, but this contradictory to our hypothesis.

ii. If $n > \eta + 1$. We have that

$$\mathrm{indC}_n(\mathrm{proj}\Lambda) \subseteq \mathrm{indC}_{\eta+1}(\mathrm{proj}\Lambda) \bigcup \mathrm{indC}_{\eta+1}(\mathrm{proj}\Lambda)[1] \bigcup \cdots \bigcup \mathrm{indC}_{\eta+1}(\mathrm{proj}\Lambda)[n-1].$$

where $\mathrm{indC}_n(\mathrm{proj}\Lambda)$ denoted the isomorphisms classes of the indecomposable complexes in $\mathrm{C}_n(\mathrm{proj}\Lambda)$ for all $n \geq 1$.

So the result follow of Lemma 2.2 and the hypothesis.

In the other hand, it is evident.

□

We give two examples in which we show how to calculate the strong global dimension of an algebra with the conditions given in the Algorithm 2.19.

**Example 3.6.** We are going to calculate the strong global dimensional of the following of finite dimensional $\mathbb{K}$-algebra. We considered $\Lambda$ the path algebra given by the following quiver

$$Q \ : \ 1 \xrightarrow{\alpha} 2 \xrightarrow{\beta} 3$$

with $\alpha\beta = 0$.

By [[15], Theorem 4.4] and [[14], Proposition 3.3] we know that $\mathrm{s.gl.dim}(\Lambda) = 2$. Now, we apply the Algorithm 2.19. First, we build the Auslander-Reiten quiver of $\mathrm{C}_2(\mathrm{proj}\Lambda)$:

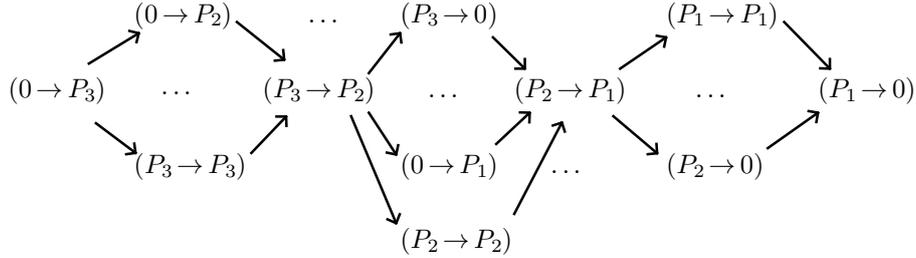

We noted that the complex $X := P_3 \longrightarrow P_2$ in $\mathrm{C}_2(\mathrm{proj}\Lambda)$ satisfies that the first and the last cells are non-zero. Therefore, we build the Auslander-Reiten quiver of $\mathrm{C}_3(\mathrm{proj}\Lambda)$:

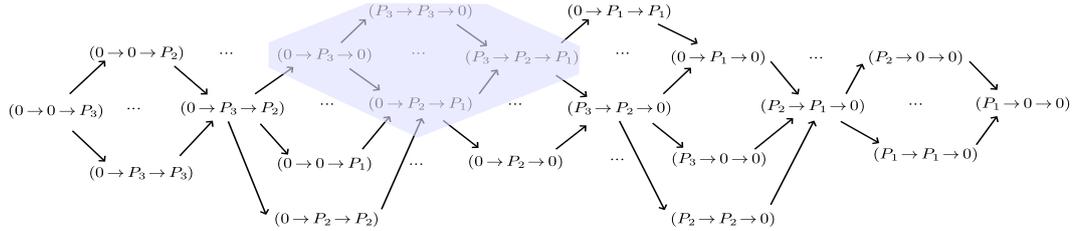

Now, the complex $X := P_3 \longrightarrow P_2 \longrightarrow P_1$ in $\mathrm{C}_3(\mathrm{proj}\Lambda)$ satisfies that the first and the last cells are non-zero, then we build the Auslander-Reiten quiver of $\mathrm{C}_4(\mathrm{proj}\Lambda)$

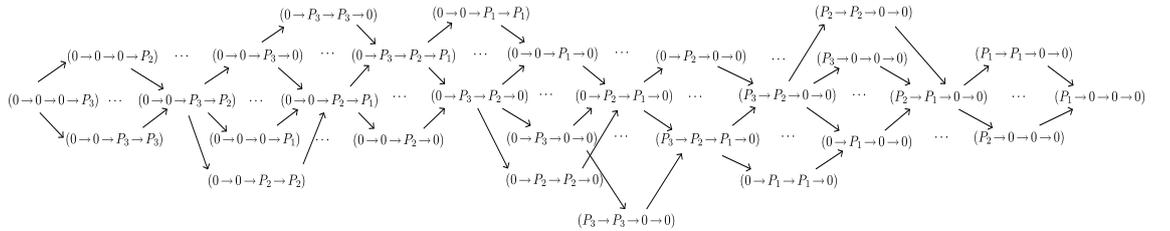

In this step all complexes in $\mathrm{C}_4(\mathrm{proj}\Lambda)$ satisfies that the first or the last cells are zero, then by Theorem 2.18 we have $\mathrm{s.gl.dim}(\Lambda) = 2$.

**Remark 3.7.** In the example above

i. Every irreducible morphism and every almost split secuence of $C_4(\text{proj}\Lambda)$ is an irreducible morphism and an almost split sequence of $C_3(\text{proj}\Lambda)$ and vice versa, that is, up to aplying abrupt truncation or inclusion respectively. This the result of the Theorem 2.16.

ii. We note that the shaded almost split sequence in the above Auslander-Reiten quiver for $C_3(\text{proj}\Lambda)$ is not an almost split sequence in the Auslander-Reiten quiver in $C_2(\text{proj}\Lambda)$, then we can say that there is one almost split secuence in $C_3(\text{proj}\Lambda)$ such that this is not an almost split secuence in $C_2(\text{proj}\Lambda)$. This fact can be generalized. There is one almost split sequence in $C_{\eta+1}(\text{proj}\Lambda)$ such that is not an almost split secuence in $C_\eta(\text{proj}\Lambda)$. Since $\eta = \text{s.gl.dim}(\Lambda)$ then there is $X \in K^b(\text{proj}\Lambda)$ an indecomposable complex non-zero such that $\ell(X) = \eta$, by the Remark 2.1 iii., $X \in C_{\eta+1}(\text{proj}\Lambda)$ and $X \notin C_\eta(\text{proj}\Lambda)$, again, by the Remark 2.1 i., $X$ is not $\mathcal{E}_{\eta+1}$-indecomposable projective, then by [[4], Theorem 8.2] there is only an almost split sequence in $C_{\eta+1}(\text{proj}\Lambda)$ such that end in $X$

$$A \longrightarrow B \xrightarrow{f} X.$$

Since $X \notin C_\eta(\text{proj}\Lambda)$, then the above almost split sequence is not in $C_\eta(\text{proj}\Lambda)$.

**Example 3.8.** In this example we calculate the strong global dimension of a finite dimensional $\mathbb{K}$-algebra using two tecniques. This example also shows that the strong global dimension is not the same as the global dimension. In [[5], Example 4.2.1], Y. Calderón-Henao calculated the s.gl.dim$(\Lambda)$ and gl.dim$(\Lambda)$ for $\Lambda$ the path algebra given by the following quiver

$$Q \;:\; 1 \xrightarrow{\alpha} 2 \xrightarrow{\beta} 3 \xrightarrow{\gamma} 4 \xrightarrow{\delta} 5 \xrightarrow{\eta} 6$$

with $\alpha\beta = \beta\gamma = \delta\eta = 0$.

That is,

$$\text{gl.dim}(\Lambda) = 3 \quad \text{and} \quad \text{s.gl.dim}(\Lambda) = 4.$$

Actually, Y. Calderón-Henao found the indecomposable complex with the longest length, which is

$$Z \colon P_6 \longrightarrow P_5 \longrightarrow P_3 \longrightarrow P_2 \longrightarrow P_1.$$

We observe that this tecnique is also used in [[9], Section 4].

Now, we can calculate the strong global dimension using the Algorithm 2.19. The Auslander-Reiten quivers of $C_n(\text{proj}\Lambda)$ for this algebra are extremely big, for this reason we do not put them here. Nevertheless we observe that the complex $Z$ is an indecomposable complex in $C_5(\text{proj}\Lambda)$ where the first and the last cells are non-zero. So we build the Auslander-Reiten quiver of $C_6(\text{proj}\Lambda)$ and every indecomposable complex in $C_6(\text{proj}\Lambda)$ the first or the last cell is zero, so by Theorem 2.18 we also have that s.gl.dim$(\Lambda) = 4$.

### Acknowledgments

The first author was supported by Beca Doctorado Nacional Colciencias (Convocatoria 727 de 2015). Also, this research was supported by CODI (Universidad de Antioquia, UdeA).

Instituto de Matemáticas, Universidad de Antioquia, Calle 67 No. 53-108, Medellín, Colombia. Instituto Tecnológico metropolitano-ITM
*Email address:* `yohny.calderon@udea.edu.co`
*Email address:* yohnycalderon8670@correo.itm.edu.co

Instituto de Matemáticas, Universidad de Antioquia, Calle 67 No. 53-108, Medellín, Colombia.
*Email address:* `cristianf.gallego@udea.edu.co`

Instituto de Matemáticas, Universidad de Antioquia, Calle 67 No. 53-108, Medellín, Colombia.
*Email address:* `hernan.giraldo@udea.edu.co`